\documentclass[12pt]{amsart}
\usepackage[left=2.5cm,right=2.5cm, top=2.5cm, bottom=2.5cm]{geometry}
\usepackage{amsmath}
\usepackage{amsfonts}
\usepackage{amssymb}
\usepackage{url}
\usepackage{color}
\usepackage{amssymb}
\usepackage{enumerate}
\usepackage{graphics}
\usepackage{graphicx}

\newtheorem{theorem}{Theorem}
\newtheorem{lemma}[theorem]{Lemma}

\newtheorem{proposition}[theorem]{Proposition}

\newtheorem{corollary}[theorem]{Corollary}

\newcommand{\Z}{{\mathbb Z}}

\newcommand{\PP}{{\mathbb P}}
\newcommand{\R}{\mathbb R}

\DeclareMathOperator{\conv}{conv}

\def\ia{{\rm I_A}}

\def\UGB{{\rm UGB}}
\def\Graver{{\mathcal G}}

\def\Orthant_j{{\mathcal O}_{j}}

\def\veb{{\ve b}}

\def\vece{{\ve e}}
\def\veg{{\ve g}}

\def\veu{{\ve u}}
\def\vev{{\ve v}}

\def\vex{{\ve x}}

\usepackage{enumerate}
\def\ve#1{\mathchoice{\mbox{\boldmath$\displaystyle\bf#1$}}
{\mbox{\boldmath$\textstyle\bf#1$}}
{\mbox{\boldmath$\scriptstyle\bf#1$}}
{\mbox{\boldmath$\scriptscriptstyle\bf#1$}}}

\newcommand{\boproof}{\textbf{Proof.} }
\newcommand{\eoproof}{\hspace*{\fill} $\square$ \vspace{5pt}}

\begin{document}
\title[Equality of Graver and Universal Gr\"obner bases]{Equality of Graver bases and universal Gr\"obner bases of colored partition identities}
\author{Tristram Bogart}
\address{Tristram Bogart, Queen's University, Kingston, Canada}
\author{Raymond Hemmecke}
\address{Raymond Hemmecke, Technische Universit\"at Munich, Germany}
\author{Sonja Petrovi\'c}
\address{Sonja Petrovi\'c, University of Illinois, Chicago, USA}
\date{}

\begin{abstract}
Associated to any vector configuration $A$ is a toric ideal encoded by vectors in the kernel of $A$. 
Each toric ideal has two special generating sets: the universal Gr\"obner basis and the Graver basis.  While the former is generally a proper subset of the latter, there are cases for which the two sets coincide. The most prominent examples among them are toric ideals of unimodular matrices. 
Equality of universal Gr\"obner basis and Graver basis is a combinatorial property of the toric ideal (or, of the defining matrix),  providing interesting information about ideals of higher Lawrence liftings of a matrix. Nonetheless, a general classification of all matrices for which both sets agree is far from known. We contribute to this task by identifying all cases with equality within two families of matrices; namely, those defining rational normal scrolls and those encoding homogeneous primitive colored partition identities.
\end{abstract}

\thanks{This project was made possible, in part, by the first and third authors' participation in the Mathematical Research Communities 2008.}

\maketitle

\vspace{-1cm}
\section{Introduction}

A vector configuration $A\in\Z^{d\times n}$ represents a toric ideal $\ia\subseteq k[x_1,\dots,x_n]$ in the following way: 
For a nonnegative vector $\veu^+$,  let $\vex^{\veu^+}$ denote the monomial $x_1^{\veu^+_1} \cdots x_n^{\veu^+_n}$. 
A vector $\veu$ in the lattice $\ker A$ corresponds to a binomial $\vex^{\veu^+} - \vex^{\veu^-}$ after writing  $\veu=\veu^+-\veu^-$ such that both $\veu^+$ and $\veu^-$ have only nonnegative coordinates and disjoint support. The toric ideal $\ia$ associated to $A$ is the set of all binomials arising from the lattice $\ker A$; that is, 
$\ia = \langle \vex^{\veu^+} - \vex^{\veu^-} \, : \, A\veu=\ve 0 \rangle . $

The ideal $\ia$ has finitely many reduced Gr\"obner bases (e.g., see \cite[Chp1]{gbcp}).  Denote by $\UGB(A)$ their union, the (minimal) \emph{universal Gr\"obner basis} of $\ia$.
In general, $\UGB(A)$ is properly contained in the \emph{Graver basis} $\Graver(A)$, consisting of all binomials $\vex^{\veu^+} -
\vex^{\veu^-} \in \ia$ for which there is no other binomial $\vex^{\vev^+} - \vex^{\vev^-} \in \ia$ such that $\vex^{\vev^+}$ divides $\vex^{\veu^+}$ and $\vex^{\vev^-}$ divides $\vex^{\veu^-}$. Such binomials are called \emph{primitive}.
Toric ideals, their generating sets and Gr\"obner bases play a prominent role in many different areas, such as algebraic geometry, commutative algebra, graph theory, integer programming, and algebraic statistics \cite{fulton}, \cite{villarreal}, \cite{gbcp}, \cite{algStatsBook}.

Universal Gr\"obner basis and Graver basis rarely coincide. It is known that they agree for some special toric ideals, including toric ideals of unimodular matrices \cite{St2,Sturmfels+Thomas:97}.  However, a general classification of all matrices for which both sets agree is not known.
We classify such matrices within two infinite families of interest, one of which defines a classical family of projective varieties. 

 Equality of the two bases provides information about higher Lawrence liftings of $A$, as discussed in \cite{HeNa}. To be more specific, for any $N\in\Z_{>0}$, consider the $N$-fold matrix
\[
[A,B]^{(N)} :=
\left(
\begin{array}{cccc}
B & B & \cdots & B \\
A & 0 &   & 0 \\
0 & A &   & 0 \\
  &   & \ddots &   \\
0 & 0 &   & A
\end{array}
\right)
\]
associated to integer matrices $A$ and $B$ of suitable dimensions. For $\veu=\left(\veu^{(1)},\ldots,\veu^{(N)}\right)\in\ker\left([A,B]^{(N)}\right)$ we call $\left|i:\veu^{(i)}\neq\ve 0\right|$ its \emph{type}. A surprising property of $N$-fold matrices is that, for all $N$, the maximum type of an element in $\Graver([A,B]^{(N)})$ is bounded by some number $g(A,B)$ \emph{not depending} on $N$. Analogously to this so-called \emph{Graver complexity} $g(A,B)$, one can define the (universal) \emph{Gr\"obner complexity} $u(A,B)$ as the maximum type of an element in $\UGB([A,B]^{(N)})$ over all $N$. As $\UGB([A,B]^{(N)})\subseteq\Graver([A,B]^{(N)})$, we have $u(A,B)\leq g(A,B)$. The main result of \cite{HeNa} is that if $\UGB(A)=\Graver(A)$, then $u(A,B)=g(A,B)$ for \emph{all} integer matrices $B$ of suitable dimensions. Hence, we have the following sequence of implications:
\[
A \text{ is unimodular } \Rightarrow \UGB(A)=\Graver(A) \Rightarrow u(A,B)=g(A,B) \text{ for all suitable matrices } B.
\]
It is known that the converse of each of these two implications is false. For example, we present infinitely many non-unimodular matrices $A$ with $\UGB(A)=\Graver(A)$. 
That the converse of the second implication does not hold has been shown 
in \cite{HeNa} for matrices of the form $A=\left(\begin{smallmatrix}1&1&1&1\\0&a&b&a+b\\\end{smallmatrix}\right)$, for which  $u(A,I_4)=g(A,I_4)$, while $\UGB(A)\subsetneq\Graver(A)$.

Let us define the two families we study. Note that the structure of the matrices resembles $N$-fold, a fact which we will exploit in our proofs. 
Given a partition  $n_1\geq n_2\geq \cdots \geq n_c$ of a positive integer $n$, define
\[
	A_{S(n_1-1,\dots,n_c-1)}  :=
	\left[
      \begin{array}{ccccccccccc}
    	1      & 2       & \cdots & n_1    & 1 & \cdots & n_2 & \cdots & 1 & \cdots & n_c\\
    	1      & 1       & \cdots & 1      & 0 & \cdots &  0   & \cdots & 0 & \cdots & 0 \\
    	0      & 0       & \cdots & 0      & 1 & \cdots & 1   & \cdots & 0 & \cdots & 0 \\
    	\vdots & \vdots  &        & \vdots &   &        &     & \ddots &   &        & \vdots \\
    	0      & 0       & \cdots & 0      & 0 & \cdots & 0   & \cdots & 1 & \cdots & 1 \\
	  \end{array}
    \right]
\]
and
\[
	A_{H(n_1,\dots,n_c)}  :=
	\left[
      \begin{array}{ccccccccccc}
    	1      & 1       & \cdots &   1    & 1 & \cdots & 1   & \cdots & 1 & \cdots & 1 \\
    	1      & 2       & \cdots & n_1    & 0 & \cdots & 0   & \cdots & 0 & \cdots & 0 \\
    	0      & 0       & \cdots & 0      & 1 & \cdots & n_2 & \cdots & 0 & \cdots & 0 \\
    	\vdots & \vdots  &        & \vdots &   &        &     & \ddots &   &        & \vdots \\
    	0      & 0       & \cdots & 0      & 0 & \cdots & 0   & \cdots & 1 & \cdots & n_c \\
	  \end{array}
    \right] .
\]
The toric ideal associated to the first matrix is the defining ideal of the $c$-dimensional \emph{rational normal scroll} $S := S(n_1-1, \ldots, n_c-1)$ in $\mathbb P^{n-1}$; see Lemma 2.1 in \cite{Pet}. In fact, the scroll $S$ is the projective variety whose defining  ideal  is generated by the $2
\times 2$ minors of the matrix
$M = \left[ M_{n_1} |  \cdots | M_{n_c} \right]$ of indeterminates, where
\begin{align*}
	 M_{n_j} = \begin{bmatrix} x_{j,1} & \cdots & x_{j,n_j-1} \\
	x_{j,2} & \cdots & x_{j,n_j} \end{bmatrix}.
\end{align*}
Eisenbud and Harris~\cite{EiHa} survey the geometry of these scrolls. In particular, they belong to the family of nondegenerate projective varieties of minimum possible degree:
just one more than the codimension.
The case $c=1$ represents the \emph{rational normal curve} in $\PP^{n-1}$.
An important combinatorial feature of this family is that the
binomials in the ideal of any rational normal scroll encode
color-homogeneous colored partition identities, as defined in
\cite{Pet}. For a precise definition, see Section \ref{sec:partitions+proof}.

The second family of matrices $A_{H(n_1,\dots,n_c)}$ represents toric ideals whose binomials encode homogeneous (but not necessarily color-homogeneous) colored partition identities.
Here, in contrast to $A_{S(n_1-1,\dots,n_c-1)}$, the positions of the
incidence vectors of the parts of the partition and the vectors $(1,2,\ldots,n_i)$ are interchanged. This effectively removes the requirement for homogeneity of the ideal with respect to the rows of the scroll matrix containing the incidence vectors.

Within these two families of matrices, we classify those for which the universal Gr\"obner and Graver bases coincide.
To state the classification, define a \emph{dominance} partial order on partitions:  $(n_1,\ldots, n_c) \preceq (m_1, \ldots, m_d)$ if $c \leq d$ and $n_j \leq m_j$ for $j=1,\ldots,c$. Naturally, this induces a partial order on $S(n_1,\ldots,n_c)$ and on $H(n_1,\dots,n_c)$: we say that $S(m_1,\ldots,m_d)$ \emph{dominates} $S(n_1,\ldots,n_c)$ if $(n_1,\ldots, n_c) \preceq (m_1, \ldots, m_d)$. Analogously, we say that $H(m_1,\ldots,m_d)$ \emph{dominates} $H(n_1,\ldots,n_c)$ if $(n_1,\ldots, n_c) \preceq (m_1, \ldots, m_d)$.

Our main result is the following.

\begin{theorem} \label{thm:minimalexamples}
The universal Gr\"obner basis of $S=S(n_1-1, \ldots, n_c-1)$ does \emph{not} equal its Graver basis if and only if $S$ dominates $S(6)$, $S(5,4)$, or $S(4,3,2)$.

The universal Gr\"obner basis of $H=H(n_1, \ldots, n_c)$ does \emph{not} equal its Graver basis if and only if $H$ dominates $H(7)$, $H(6,2)$, or $H(4,3)$.
\end{theorem}

This result has some interesting consequences. First,
it answers a question posed in \cite{HeNa}, showing that $N$-fold matrices that preserve the Gr\"obner and Graver complexities do not preserve equality of universal Gr\"obner and Graver bases. Namely, 
 there exist integer matrices $A$ and $B$ of appropriate dimensions satisfying $\UGB(A)=\Graver(A)$, and hence $u(A,B)=g(A,B)$, but that still satisfy $\UGB([A,B]^{(N)})\subsetneq\Graver([A,B]^{(N)})$ for all $N>1$. 
 One such example is given by $S(5)$, that is, $A=(1,1,1,1,1,1)$ and $B=(1,2,3,4,5,6)$.

Secondly, let us point out that $\UGB(A)=\Graver(A)$ seems to be a purely \emph{combinatorial} property of $A$. It is \emph{not} related to the sequence of {algebraic} properties
\[
A \text{ is unimodular } \Rightarrow A \text{ is compressed } \Rightarrow K\text{-algebra } K[A] \text{ is normal (for some field } K),
\]
as discussed in \cite{Hibi+Ohsugi}. Note that the first family of matrices is not compressed. 
In addition, among the matrices within the two families that satisfy $\UGB(A)=\Graver(A)$, one set has this normality property, while the other one doesn't. It is well-known that semigroup algebras $K[A]$ of rational normal scrolls are normal (this follows, for example, by Proposition 13.5 in \cite{gbcp}, since every scroll has a squarefree initial ideal). 
On the other hand, the second family of matrices does not have this property:
\begin{lemma}
$K[H]:=K[H(n_1,\ldots,n_c)]$ is normal if and only if $H$ does not dominate $H(2,2)$.
\end{lemma}
\boproof
 Let $H(2,2)=\left(\begin{smallmatrix}1&1&1&1\\1&2&0&0\\0&0&1&2\\\end{smallmatrix}\right)$. 
 Then $K[H(2,2)]$ is not normal by Lemma 6.1 in \cite{Hibi+Ohsugi}, since the binomial $\vex_1^2\vex_4-\vex_2\vex_3^2\in I_{H(2,2)}$ is an indispensible binomial. 
As we find a similar binomial (after suitable index permutation) in the toric ideal of any $H$ dominating $H(2,2)$, one implication of the claim follows.
It remains to show that if $H$ does not dominate $H(2,2)$ then $K$-algebra of $H$ is normal. But then $H$ is dominated by $H(k,1,\ldots,1)$ and hence the corresponding $K$-algebra is isomorphic to the $K$-algebra of $S(k-1)$ which is normal. \eoproof

As we will see in the following section (in particular, Corollary \ref{Corollary: Equality is inherited by submatrices}),
Theorem \ref{thm:minimalexamples} gives the classification also for certain submatrices of the matrices we consider here.
 It is not clear whether, in case both sets do not coincide, the ratios $\left| \UGB({A_S}) \right| / \left| \Graver(A_S) \right| $ and $\left| \UGB({A_H}) \right| / \left| \Graver(A_H) \right| $ tend to $0$ as $S$ and $H$ increase in the dominance ordering.
Moreover, it is an interesting open question whether the special structure of our matrices implies that there are only finitely many fundamental counterexamples to equality of universal Gr\"obner basis and Graver basis, from which all other counterexamples can be derived.

In particular, such results would provide some insight into the complexity of the \emph{Gr\"obner fan} of varieties of minimal degree.
Namely, each reduced Gr\"obner basis determines a cone in the {Gr\"obner fan} of the ideal $\ia$.
 A part of the Gr\"obner fan of rational normal curves is understood: in  \cite{ConcaDeNegriRossi}, Conca, De Negri and Rossi determine explicitly all the initial ideals of $I_{S(n)}$ that are Cohen-Macaulay.
To determine the rest of the fan for such a curve or for a general scroll remains an open problem.
One approach to this problem is to first understand an approximation of the set $\UGB(A)$.  Since the Graver basis equals the universal Gr\"obner basis for only a small set of scrolls, understanding the primitive binomials that do not belong to $\UGB(A)$ is of interest.

\section{Preliminaries}

Universal Gr\"obner bases and Graver bases have nice geometric properties that allow us to enumerate  the dominance-minimal cases where equality fails.
In this section we recall known results fundamental to our problem.

For any $\veb\in\Z^d$, the polytope $P_\veb^I := \conv(\{\veu:A\veu=\veb,\veu\in\Z^n_+\})$ is called the \emph{fiber} of $\veb$. With this, we can characterize the elements in $\UGB(A)$.

\begin{proposition} \cite{gbcp} \label{prop:groebnerfibers}
A integer vector $\veu^+ - \veu^-$ (with $\veu^+, \veu^- \geq\ve 0$) lies in $\UGB(A)$  if and only if the line segment $[\veu^+, \veu^-]$ is an edge of the fiber $P_{A \veu^+}^{I}$ and contains no lattice points other than its endpoints.
\end{proposition}

This proposition allows us to check  computationally whether
$\UGB(A)=\Graver(A)$ for any given matrix $A$ as follows: First, we
compute the Graver basis $\Graver(A)$;  then, we test the statement in
Proposition \ref{prop:groebnerfibers} for every
$\veu\in\Graver(A)$. To check this condition, we first enumerate the
set $V$ of all lattice points in the polytope $P_{A \veu^+}$ and then check whether $\veu$ is an edge of this polytope that contains no interior lattice point. This can be done by first testing whether $\veu^+$ and $\veu^-$ are vertices of $P_{A \veu^+}$ (that is, they are not convex combination of $V\setminus\{\veu^-\}$ and of $V\setminus\{\veu^+\}$, respectively). If $\veu^+$ and $\veu^-$ are vertices of $P_{A \veu^+}$, then $\veu$ is an edge with no integer point in its interior if and only if there does not exist a decomposition
\[
\veu=\veu^+-\veu^-=\sum_{\vev\in V\setminus\{\veu^+,\veu^-\}} \lambda_\vev (\vev-\veu^-)
\]
with non-negative real coefficients $\lambda_\vev$. This feasibility problem can be decided by any code that solves linear programs. We applied the commercial solver CPLEX \cite{Cplex}.

For small scrolls, the universal Gr\"obner basis can also be computed
by the software package Gfan~\cite{gfan}. However, once we reached $9$
variables, the program failed to finish this computation due to its
complexity. Even when the universal Gr\"obner basis itself is not
especially large, each binomial occurs in many different Gr\"obner
bases, each one of which must be enumerated in order to compute the
Gr\"obner fan, as Gfan does. In contrast, 4ti2 \cite{4ti2} is very
quick in calculating the Graver basis for our two families of matrices. We then use
Proposition~\ref{prop:groebnerfibers} to extract the universal Gr\"obner basis from the Graver basis set. Such an approach was necessary to make the computations feasible for reasonably large examples.

Besides this algorithmic test, Proposition \ref{prop:groebnerfibers} also allows us to show a well-known and very useful fact; namely, that certain projections preserve elements in the universal Gr\"obner and Graver bases:

\begin{corollary} \label{cor:extendUGB}
Suppose $\veu \in \ker_{\Z} A$ and $\sigma \subseteq [n]$ is such that $\veu_i = 0$ for $i \notin \sigma$.  Let $A_{\sigma}$ be the submatrix of $A$ of columns indexed by $\sigma$ and $\veu_{\sigma}$ be the projection of $\veu$ to $\R^{\sigma}$. Then:
\begin{itemize}
  \item[(a)] $\veu \in \UGB(A)$ if and only if $\veu_{\sigma} \in \UGB(A_{\sigma})$.
  \item[(b)] $\veu \in \Graver(A)$ if and only if $\veu_{\sigma} \in \Graver(A_{\sigma})$.
\end{itemize}
\end{corollary}

\boproof
To prove claim (a), observe that since the hyperplanes $x_i = 0$ ($i \notin \sigma$) do not pass through the interior of $P_{A \veu^+}$, the polytope $P_{A_{\sigma} \veu_{\sigma}^+}$ is a face of $P_{A \veu^+}$. It follows that the segment $[\veu^+, \veu^-]$ is an edge of $P_{A_{\sigma} \veu_{\sigma}^+}$ if and only if it is an edge of $P_{A \veu^+}$.

To prove claim (b), simply observe that the minimality property of
$\vex^{\vev^+} - \vex^{\vev^-} \in \ia$ is the same for both matrices
$A$ and $A_{\sigma}$, as the variables indexed by elements of $\sigma$ do not appear in either case.
\eoproof

Note that Corollary~\ref{cor:extendUGB} immediately implies that equality of the two sets is inherited by submatrices:

\begin{corollary}\label{Corollary: Equality is inherited by submatrices}
Let $A_{\sigma}$ be obtained from $A$ by first choosing the submatrix of $A$ consisting of the columns indexed by $\sigma \subseteq [n]$ and by then eliminating some or all of the redundant rows. Then $\UGB(A)=\Graver(A)$ implies $\UGB(A_{\sigma})=\Graver(A_{\sigma})$.
\end{corollary}

The dominance order on partitions provides the following simple consequence of Corollary \ref{Corollary: Equality is inherited by submatrices}  (See also Proposition 4.13 in \cite{gbcp}).

\begin{proposition}\label{Proposition: dominated scrolls inherit equality}
Suppose  $(n_1,\ldots n_c) \prec (m_1, \ldots, m_d)$. Then $\UGB\left(A_{S(m_1, \ldots, m_d)}\right)=\Graver\left(A_{S(m_1, \ldots, m_d)}\right)$ implies $\UGB\left(A_{S(n_1, \ldots, n_c)}\right)=\Graver\left(A_{S(n_1, \ldots, n_c)}\right)$. Similarly, $\UGB\left(A_{H(m_1, \ldots, m_d)}\right)=\Graver\left(A_{H(m_1, \ldots, m_d)}\right)$ implies $\UGB\left(A_{H(n_1, \ldots, n_c)}\right)=\Graver\left(A_{H(n_1, \ldots, n_c)}\right)$.
\end{proposition}

In particular, this allows us to solve our classification problem by listing only the dominance-minimal matrices for which equality of the universal Gr\"obner basis and Graver basis does \emph{not} hold.

\section{Partitions, Graver Complexity and the Proof of Theorem~\ref{thm:minimalexamples}}
\label{sec:partitions+proof}

The toric ideals in the two families we are studying have a nice combinatorial interpretation.
The binomials $x_{a_1}x_{a_2} \cdots x_{a_k} - x_{b_1}x_{b_2} \cdots x_{b_k}$ in the ideal $I_{S(n_1-1)}$ of a rational normal curve
encode \emph{homogeneous partition identities}:
$$a_1 + \cdots + a_k = b_1 + \cdots + b_k$$
where $a_1, \ldots, a_k, b_1, \ldots, b_k$ are positive
integers, not necessarily distinct \cite{DGS}. Similarly, the
binomials in the ideal $I_S$ of a rational normal scroll encode
\emph{color-homogeneous {c}olored {p}artition {i}dentities} (color-homogeneous cpi's) \cite{Pet}:
$$ a_{1,1} + \cdots + a_{1, k_1} + \cdots + a_{c, 1} + \cdots + a_{c,k_c} = b_{1,1} + \cdots + b_{1, k_1} + \cdots + b_{c, 1} + \cdots + b_{c, k_c}.$$

For example,
$${\color{red} 1_1} + {\color{red}5_1} + {\color{blue}1_2} + {\color{blue}5_2} = {\color{red}2_1} + {\color{red}6_1} + {\color{blue}2_2} + {\color{blue}2_2} $$
is a color-homogeneous cpi encoded by the binomial
\[
	\textcolor{red}{x_{1,1}x_{1,5}}\textcolor{blue}{x_{2,1}x_{2,5}}- \textcolor{red}{x_{1,2}x_{1,6}}\textcolor{blue}{x_{2,2}^2}	,
\]
while
$${\color{red} 1_1} + {\color{red}5_1} + {\color{blue}1_2} = {\color{red}3_1} + {\color{blue}1_2} + {\color{blue}3_2}$$
is a cpi that is homogeneous, but not color-homogeneous, and is encoded by the binomial
\[
	\textcolor{red}{x_{1,1}x_{1,5}}\textcolor{blue}{x_{2,1}}-\textcolor{red}{x_{1,3}}\textcolor{blue}{x_{2,1}x_{2,3}}.
\]

Among such identities, those with no proper sub-identities are again called \emph{primitive} and are encoded by the elements of $\Graver(A_S)$ \cite{Pet}.
The analogous statement is true for $\Graver(A_H)$.
This attractive characterization of the Graver bases makes it especially useful to classify which scrolls have the property that the universal Gr\"obner basis is equal to the Graver basis.

Let us now start collecting the ingredients for the proof of our main result.

\begin{lemma}\label{Lemma: Counter-examples to the three minimal cases}
The universal Gr\"obner basis and the Graver basis are not the same for the defining matrices of $S(6)$, $S(5,4)$, $S(4,3,2)$, $H(7)$, $H(6,2)$, and $H(4,3)$.
\end{lemma}

\boproof
The defining matrix of $S(6)$ is
\[
A_{S(6)}=\left(
\begin{array}{ccccccc}
1 & 2 & 3 & 4 & 5 & 6 & 7\\
1 & 1 & 1 & 1 & 1 & 1 & 1\\
\end{array}
\right).
\]
Consider the vector $\veg=(1,-1,1,-1,-1,0,1)\in\ker(A_{S(6)})$ and the three vectors $\veu_1=(1,0,0,0,2,0,0)$, $\veu_2=(0,2,0,0,0,0,1)$, $\veu_3=(0,0,1,2,0,0,0)$ in the fiber $\{\veu:A_{S(6)}\veu=A_{S(6)}\veg^+,\veu\in\Z^7_+\}$ of $\veg$. Notice that
\[
\veg=\veg^+-\veg^-=1\cdot(\veu_1-\veg^-)+1\cdot(\veu_2-\veg^-)+1\cdot(\veu_3-\veg^-)
\]
and thus $\veg$ is not an edge of $P_{\veg^+}^I.$
Consequently, $\vex^{\veg^+} - \vex^{\veg^-} \not\in\UGB(A_{S(6)})$. Yet $\veg$ represents the homogeneous primitive partition identity
$$1+3+7=2+4+5$$ and thus $\vex^{\veg^+} - \vex^{\veg^-} \in \Graver(A_{S(6)}) \setminus \UGB(A_{S(6)})$.

The defining matrix of $S(5,4)$ is
\[
A_{S(5,4)}=\left(
\begin{array}{ccccccccccc}
1 & 2 & 3 & 4 & 5 & 6 & 1 & 2 & 3 & 4 & 5 \\
1 & 1 & 1 & 1 & 1 & 0 & 0 & 0 & 0 & 0 & 0 \\
0 & 0 & 0 & 0 & 0 & 1 & 1 & 1 & 1 & 1 & 1 \\
\end{array}
\right).
\]
Consider $\veg=(1,-1,0,0,1,-1,1,-2,0,0,1)\in\ker(A_{S(5,4)})$,
$\veu_1=(0,0,0,0,2,0,2,0,0,0,0)$, and $\veu_2=(2,0,0,0,0,0,0,0,0,0,2)$. Note that $\veu_1,\veu_2 \in \{\veu:A_{S(5,4)}\veu=A_{S(5,4)}\veg^+,\veu\in\Z^{11}_+\}$.
Then
\[
\veg^+=\frac{1}{2}\cdot\veu_1+\frac{1}{2}\cdot\veu_2,
\]
and thus $\veg^+$ is not a vertex of $P_{\veg^+}^I.$
Consequently, $\vex^{\veg^+} - \vex^{\veg^-} \not\in\UGB(A_{S(5,4)})$.
Yet $\veg$ represents the color-homogeneous primitive partition identity
$${\color{red} 1_1} + {\color{red}5_1} + {\color{blue}1_2} + {\color{blue}5_2} = {\color{red}2_1} + {\color{red}6_1} + {\color{blue}2_2} + {\color{blue}2_2} $$
and thus $\vex^{\veg^+} - \vex^{\veg^-} \in \Graver(A_{S(5,4)}) \setminus \UGB(A_{S(5,4)})$.

The defining matrix of $S(4,3,2)$ is
\[
A_{S(4,3,2)}=\left(
\begin{array}{cccccccccccc}
1 & 2 & 3 & 4 & 5 & 1 & 2 & 3 & 4 & 1 & 2 & 3 \\
1 & 1 & 1 & 1 & 1 & 0 & 0 & 0 & 0 & 0 & 0 & 0 \\
0 & 0 & 0 & 0 & 0 & 1 & 1 & 1 & 1 & 0 & 0 & 0 \\
0 & 0 & 0 & 0 & 0 & 0 & 0 & 0 & 0 & 1 & 1 & 1 \\
\end{array}
\right).
\]
Consider the vectors
\begin{eqnarray*}
\veg & = &(-1,0,0,0,1,1,-1,1,-1,1,0,-1) \in \ker(A_{S(4,3,2)}),\\
\veu_1 & = & (0,0,0,0,1,0,2,0,0,1,0,0),\\
\veu_2 & = & (1,0,0,0,0,0,0,0,2,1,0,0),\\
\veu_3 & = & (0,0,0,0,1,2,0,0,0,0,0,1), \text{and}\\
\veu_4 & = & (1,0,0,0,0,0,0,2,0,0,0,1).
\end{eqnarray*}
Note that $\veu_1,\veu_2,\veu_3,\veu_4\in\{\veu:A_{S(4,3,2)}\veu=A_{S(4,3,2)}\veg^+,\veu\in\Z^{12}_+\}$. The decomposition
\[
\veg=\veg^+-\veg^-=\frac{1}{2}\cdot(\veu_1-\veg^-)+\frac{1}{2}\cdot(\veu_2-\veg^-)+\frac{1}{2}\cdot(\veu_3-\veg^-)+\frac{1}{2}\cdot(\veu_4-\veg^-),
\]
implies that $\veg$ is not an edge of $P_{\veg^+}^I$
and, consequently, $\vex^{\veg^+} - \vex^{\veg^-} \not\in\UGB(A_{S(4,3,2)})$.
Yet $\veg$ represents the color-homogeneous primitive partition identity
$${\color{red}5_1}+{\color{blue}1_2}+{\color{blue}3_2}+{\color{green}1_3} = {\color{red}1_1}+{\color{blue}2_2}+{\color{blue}4_2}+{\color{green}3_3}$$
and thus $\vex^{\veg^+} - \vex^{\veg^-} \in \Graver(A_{S(4,3,2)}) \setminus \UGB(A_{S(4,3,2)})$.

The defining matrix of $H(7)$ is the same as of $S(6)$ (up to a swap of rows). Thus the toric ideal corresponding to $H(7)$ is the same as for $S(6)$ and the counterexample for $S(6)$ applies here, too.

The defining matrix of $H(6,2)$ is
\[
A_{H(6,2)}=\left(
\begin{array}{cccccccc}
1 & 1 & 1 & 1 & 1 & 1 & 1 & 1 \\
1 & 2 & 3 & 4 & 5 & 6 & 0 & 0 \\
0 & 0 & 0 & 0 & 0 & 0 & 1 & 2 \\
\end{array}
\right).
\]
Consider $\veg=(-1,1,-1,-1,0,1,2,-1)\in\ker(A_{H(6,2)})$,
$\veu_1=(0,0,0,2,2,0,0,0)$, $\veu_2=(2,0,0,0,0,1,0,1)$, and $\veu_3=(0,1,2,0,0,0,1,0,0)$. Note that $\veu_1,\veu_2,\veu_3 \in \{\veu:A_{H(6,2)}\veu=A_{H(6,2)}\veg^+,\veu\in\Z^{8}_+\}$.
Again,  $\veg$ is not an edge of $P_{\veg^+}^I$, since
\[
\veg=\veg^+-\veg^-=1\cdot(\veu_1-\veg^-)+1\cdot(\veu_2-\veg^-)+1\cdot(\veu_3-\veg^-).
\]
Consequently, $\vex^{\veg^+} - \vex^{\veg^-} \not\in\UGB(A_{H(6,2)})$.
Yet $\veg$ represents the homogeneous primitive colored partition identity
$${\color{red} 2_1} + {\color{red} 6_1} + {\color{blue}1_2} + {\color{blue}1_2} = {\color{red}1_1} + {\color{red}3_1} + {\color{red}4_1} + {\color{blue}2_2}$$
and thus $\vex^{\veg^+} - \vex^{\veg^-} \in \Graver(A_{H(6,2)}) \setminus \UGB(A_{H(6,2)})$.

The defining matrix of $H(4,3)$ is
\[
A_{H(4,3)}=\left(
\begin{array}{ccccccc}
1 & 1 & 1 & 1 & 1 & 1 & 1 \\
1 & 2 & 3 & 4 & 0 & 0 & 0 \\
0 & 0 & 0 & 0 & 1 & 2 & 3 \\
\end{array}
\right).
\]
Consider $\veg=(1,2,1,-2,-3,0,1)\in\ker(A_{H(4,3)})$,
$\veu_1=(2,0,2,0,0,0,1)$, and $\veu_2=(0,4,0,0,0,0,1)$. Note that $\veu_1,\veu_2 \in \{\veu:A_{H(4,3)}\veu=A_{H(4,3)}\veg^+,\veu\in\Z^{7}_+\}$.
This time, $\veg^+$ is not a vertex of $P_{\veg^+}^I$ since
\[
\veg^+=\frac{1}{2}\cdot\veu_1+\frac{1}{2}\cdot\veu_2.
\]
On the other hand,  $\veg$ represents the homogeneous primitive colored partition identity
$${\color{red} 1_1} + {\color{red}2_1} + {\color{red}2_1} + {\color{red}3_1} + {\color{blue}3_2} = {\color{red}4_1} + {\color{red}4_1} + {\color{blue}1_2} + {\color{blue}1_2} + {\color{blue}1_2}, $$
and thus $\vex^{\veg^+} - \vex^{\veg^-} \in \Graver(A_{H(4,3)}) \setminus \UGB(A_{H(4,3)})$.
\eoproof

Consider the families $S_{c,m-1}:=S(m-1,m-1,\ldots,m-1)$ and $H_{c,m}:=S(m,m,\ldots,m)$ with $c$ components $m-1$ and $m$, respectively.
The corresponding matrices $A_{c,m-1}:=A_{S(m-1,\ldots,m-1)}$ and $B_{c,m}:=A_{H(m,\ldots,m)}$ have a special structure in this case: they are $c$-fold matrices of the form $A_{c,m-1} := [C,D]^{(c)}$ and $B_{c,m} :=[D,C]^{(c)}$ with $C=(1,1,\ldots,1)$ and $D=(1,2,\ldots,m)$.

The special structure of these matrices allows us to apply known results about generalized higher Lawrence liftings \cite{Hosten+Sullivant} that the $c$-fold matrix represents. In particular, this structure implies that, as $c$ grows, the Graver bases of $A_{c,m}$ eventually stabilize (see \cite{Hosten+Sullivant}), in the sense that the type $\left|i:\veu^{(i)}\neq\ve 0\right|$ of a vector $\veu=\left(\veu^{(1)},\ldots,\veu^{(c)}\right)\in\Graver\left(A_{c,m}\right)$ is bounded by a constant $g(C,D)$, the so-called \emph{Graver complexity} of $C$ and $D$. A similar bound $g(D,C)$ exists for the type of the Graver basis elements of $B_{c,m}$. Note, however, that generally $g(C,D)\neq g(D,C)$.

\begin{lemma} \label{lem:complexity}
For any fixed $m$,  $C=(1,1,\ldots,1)$, and $D=(1,2,\ldots,m)$,
the Graver complexities of the $c$-fold matrices satisfy
 $g(C,D)=2m-3$ and $g(D,C)\leq 4m-7$.
\end{lemma}

\boproof
The Graver complexity $g(C,D)$ can be computed via the algorithm presented in \cite{Hosten+Sullivant}:
\[
g(C,D)=\max\{\|\vev\|_1:\vev\in\Graver(D\cdot\Graver(C))\}.
\]
The Graver basis of the $1\times m$ matrix $C=(1,1,\ldots,1)$ consists of all vectors $\vece_i-\vece_j$, $1\leq i\neq j\leq m$. Multiplying these elements by the $1\times m$ matrix $D=(1,2,3,\ldots,m)$, we conclude that the desired Graver complexity equals the maximum $1$-norm among the Graver basis elements of the matrix $(1,2,3,\ldots,m-1)$, which is known to be $2(m-1)-1=2m-3$ \cite{DGS}.

Similarly, we can compute $g(D,C)$. The Graver basis elements of $D=(1,2,3,\ldots,m)$ have a maximum $1$-norm of $2m-1$. However, as no element in $\Graver(D)$ has only nonnegative entries, we have that $|C\veg|\leq 2m-3$. Hence, $g(D,C)$ is bounded from above by the maximum $1$-norm among the Graver basis elements of the matrix $(1,2,3,\ldots,2m-3)$, which is $2(2m-3)-1=4m-7$.
\eoproof

This lemma can be exploited computationally as it bounds the sizes $c$ of $A_{c,m}$ and $B_{c,m}$ for which equality of the universal Gr\"obner basis and the Graver basis have to be checked. We start with two simple cases.

\begin{lemma}\label{Lemma: Equality holds for A_{c,3}}
The universal Gr\"obner basis equals the Graver basis for $A_{c,3}$ and for $B_{c,3}$ for any $c\in\Z_{>0}$.
\end{lemma}

\boproof
We first verify computationally that the universal Gr\"obner basis equals the Graver basis for $S(3,3,3,3,3)$ and for $H(3,3,3,3,3)$. Now suppose $\veg \in \Graver(A_{c,3})$ for some $c > 5$. By Lemma~\ref{lem:complexity}, the type of any Graver basis element of $A_{c,3}$ is at most $2\cdot 4-3=5$, so $\veg$ represents a color-homogeneous primitive partition identity with at most five colors. Define $\veg'$ by restricting $\veg$ to the $20$ coordinates that represent these five colors (eliminating only zeros.)
Then $\veg' \in \Graver(A_{5,3}) = \UGB(A_{5,3})$. Thus, by applying Corollary~\ref{cor:extendUGB}, we see that $\veg \in \UGB(A_{c,3})$.

Similar arguments apply to $B_{c,3}$. Here, the type of any Graver
basis element is bounded by $4\cdot 3-7=5$, and the result follows
from the equality of the universal Gr\"obner basis and Graver basis for $B_{5,3}$.
\eoproof

\begin{lemma}\label{Lemma: Equality holds for 53k1, 5k2, 44k1}
For the matrices of $S(5)$, $S(5,3,1,\ldots,1)$, $S(5,2,\ldots,2)$, $S(4,4,1,\ldots,1)$, $H(6)$, and $H(5,2,\ldots,2)$, the universal Gr\"obner basis equals the Graver basis for any number of $1$'s and $2$'s, respectively.
\end{lemma}

\boproof
The cases $S(5)$ and $H(6)$ correspond to the same toric ideal and equality can be verified computationally using Proposition \ref{prop:groebnerfibers}. Moreover, we verify computationally that the universal Gr\"obner basis equals the Graver basis for the following cases: $S(5,3,1,1,1,1,1,1,1)$, $S(5,2,2,2,2,2,2,2,2)$, $S(4,4,1,1,1,1,1)$, and $H(5,2,\ldots,2)$ (with $12$ components $2$).
To see that the result of the lemma now follows for any number of $1$'s and $2$'s, let us give the arguments for $H(5,2,\ldots,2)$. The other cases can be handled analogously.

Let $H=H(5,2,\ldots,2)$ with $k$ components $2$. If $k\leq 12$, the result follows from the equality for $k=12$ and Proposition \ref{Proposition: dominated scrolls inherit equality}. Let $k>12$ and let $\veg\in\Graver(A_{H(5,2,\ldots,2)})$. First observe that $A_{H(5,2,\ldots,2)}$ can be obtained from $A_{H(5,5,\ldots,5)}$ by deleting certain columns and zero rows thereafter. Thus, $\veg$ can be lifted (by only adding components $0$) to $\veg'\in\Graver(A_{H(5,5,\ldots,5)})$ by Corollary \ref{cor:extendUGB}. The type of any Graver basis element of $A_{H(5,5,\ldots,5)}$ is bounded by $4\cdot 5-7=13$. This implies that the type of $\veg'$ is at most $13$. Therefore, it can be projected (by only removing zero components) to a Graver basis element $\veg''\in\Graver(A_{H(5,2,\ldots,2)})$ for $12$ components $2$. As we have verified computationally that $\UGB(A_{H(5,2,\ldots,2)})=\Graver(A_{H(5,2,\ldots,2)})$ in this case, we conclude that $\veg''\in\UGB(A_{H(5,2,\ldots,2)})$ for $12$ components $2$. Applying Corollary \ref{cor:extendUGB} twice, we conclude that $\veg\in\UGB(A_{H(5,2,\ldots,2)})$ with $k>12$ components $2$, and thus, $\UGB(A_H)=\Graver(A_H)$ as claimed. \eoproof

Let us finally prove our main theorem.

{\bf Proof of Theorem \ref{thm:minimalexamples}.}

Suppose that $S=S(n_1-1, \ldots, n_c-1)$ is any scroll. If $S$
dominates $S(6)$, $S(5,4)$, or $S(4,3,2)$, then the universal Gr\"obner basis and Graver basis do not agree by Proposition \ref{Proposition: dominated scrolls inherit equality} and Lemma \ref{Lemma: Counter-examples to the three minimal cases}. Thus, let us assume that $S$ does not dominate $S(6)$, $S(5,4)$, or $S(4,3,2)$. Then one of the following four cases applies.

Case 1: $n_1 \leq 4$.\\
\mbox{}\hspace{1cm} Then $S \preceq S_{c,3}$ for some $c$ and thus equality holds by Lemma \ref{Lemma: Equality holds for A_{c,3}}.

Case 2: $n_1 = n_2 = 5$.\\
\mbox{}\hspace{1cm} Then $n_3 \leq 2$ to avoid dominating $S(4,3,2)$. Thus, we have $S=S(4,4,1, \ldots, 1)$.

Case 3: $n_1 = 5 \textup{ or } 6, n_2 \leq 3$.\\
\mbox{}\hspace{1cm} Then $S \preceq S(5,2,\ldots,2)$.

Case 4: $n_1=5 \textup{ or } 6, n_2 = 4$.\\
\mbox{}\hspace{1cm}Then $n_3 \leq 2$ to avoid dominating $S(4,3,2)$. Thus, we have $S \preceq S(5,3,1,\ldots, 1)$.

By Proposition \ref{Proposition: dominated scrolls inherit equality} and Lemma \ref{Lemma: Equality holds for 53k1, 5k2, 44k1}, the universal Gr\"obner basis and Graver basis coincide for the scrolls  $S$ in Cases 2, 3, and 4.

Suppose now that $H=H(n_1, \ldots, n_c)$. If $H$ dominates $H(7)$, $H(6,2)$, or $H(4,3)$, then universal Gr\"obner basis and Graver basis do not agree by Proposition \ref{Proposition: dominated scrolls inherit equality} and Lemma \ref{Lemma: Counter-examples to the three minimal cases}. Thus, let us assume that $H$ does not dominate $H(7)$, $H(6,2)$, or $H(4,3)$. Then one of the following three cases applies.

Case 1: $n_1 \leq 3$.\\
\mbox{}\hspace{1cm} Then $H \preceq H_{c,3}$ for some $c$ and thus equality holds by Lemma \ref{Lemma: Equality holds for A_{c,3}}.

Case 2: $n_1 = 4 \textup{ or } 5$.\\
\mbox{}\hspace{1cm} Then $n_2 \leq 2$ to avoid dominating $H(4,3)$. Then $H \preceq H(5,2,\ldots,2)$.

Case 3: $n_1 = 6$.\\
\mbox{}\hspace{1cm}Then $n_2 \leq 1$ to avoid dominating $H(6,2)$. But then $H$ has the same toric ideal as $H(6)$.

By Proposition \ref{Proposition: dominated scrolls inherit equality} and Lemma \ref{Lemma: Equality holds for 53k1, 5k2, 44k1}, the universal Gr\"obner basis and Graver basis coincide for $H$ in Cases 2 and  3.
\eoproof


\begin{thebibliography}{99}

\bibitem{ConcaDeNegriRossi}
A.~Conca, E.~De Negri, and M. E. Rossi.
Contracted ideals and the Gr\"obner fan of the rational normal curve.
Algebra Number Theory 1 (2007), no. 3, 239--268.   

\bibitem{Cplex} CPLEX Callable Library 9.1.3, ILOG (2005)

\bibitem{DGS}
P.~Diaconis, R.~Graham, and B.~Sturmfels.
Primitive partition identities.
In {\sl Combinatorics, Paul Erdos is Eighty}, eds. D.~Mikl\'os, V.~T.~S\'os, T.~Szonyi,
Janos Bolyai Mathematical Society, Budapest, Hungary, 1996, 173--192.

\bibitem{algStatsBook}
M.~Drton, B.~Sturmfels and S.~Sullivant.
 \emph{Lectures on algebraic
  statistics}, Oberwolfach Seminars \textbf{39}, Birkh\"auser  (2009)

\bibitem{EiHa}
D.~Eisenbud and J.~Harris.
On varieties of minimal degree (a centennial account).
In {\sl Algebraic Geometry, Bowdoin, 1985},
Proceedings of Symposia in Pure Mathematics,
\textbf{46(1)} (1987), 3--13.

\bibitem{fulton}
W.~Fulton. Introduction to toric varieties.
Ann. of Math. Stud. 131, Princeton University Press, Princeton, NJ, 1993.

\bibitem{4ti2} 4ti2~team.
4ti2---A software package for algebraic, geometric and combinatorial problems on linear spaces. {A}vailable at www.4ti2.de.

\bibitem{HeNa}
R.~Hemmecke and K.~Nairn.
On the {G}r\"obner complexity of matrices.
{\sl Journal of Pure and Applied Algebra} \textbf{213} (2009), 1558--1563.

\bibitem{Hibi+Ohsugi}
T.~Hibi and H.~Ohsugi.
Toric ideals arising from contingency tables.
In: Commutative Algebra and Combinatorics.
In: {\sl Ramanujan Mathematical Society Lecture Note Series} \textbf{4} (2006), 87--111.

\bibitem{Hosten+Sullivant}
S.~Ho\c sten and S.~Sullivant.
Finiteness theorems for Markov bases of hierarchical models.
{\sl Journal of Combinatorial Theory}, Series A \textbf{114(2)} (2007), 311--321.

\bibitem{gfan}
A.~N.~Jensen.
Gfan, a software system for {G}r{\"o}bner fans and tropical varieties. {A}vailable at \url{http://www.math.tu-berlin.de/~jensen/software/gfan/gfan.html}

\bibitem{Pet}
S.~Petrovi{\'c}.
\newblock On the universal {G}r{\"o}bner bases of varieties of minimal degree.
\newblock {\sl Math. Res. Lett.},
\textbf{15(6)}(2008), 1211--1221.
\newblock arXiv:0711.2714.

\bibitem{St2}
B.~Sturmfels.
\newblock Asymptotic analysis of toric ideals.
\newblock {\sl Memoirs of the faculty of science}, Kyushu University, Ser. A \textbf{46} (1992), 217--228.

\bibitem{gbcp}
B.~Sturmfels.
\newblock {\em Gr\"obner Bases and Convex Polytopes}, volume~8 of {\em
  University Lecture Series}.
\newblock American Mathematical Society, Providence, RI, 1996.

\bibitem{Sturmfels+Thomas:97}
B.~Sturmfels and R.~R.~Thomas.
\newblock Variation of cost functions in integer programming.
\newblock {\sl Mathematical Programming} \textbf{77} (1997), 357--387.

\bibitem{villarreal}
R.~Villarreal.
Monomial Algebras, Monographs and Textbooks in Pure and Applied Mathematics 238, Marcel Dekker, Inc., New York, 2001.

\end{thebibliography}
\end{document}